\documentclass[final,3p]{elsarticle}

\makeatletter
\def\ps@pprintTitle{%
 \let\@oddhead\@empty
 \let\@evenhead\@empty
 \def\@oddfoot{\centerline{\thepage}}%
 \let\@evenfoot\@oddfoot}
\makeatother

\usepackage{graphicx}
\usepackage[colorlinks=true,bookmarks,bookmarksopen,bookmarksdepth=2]{hyperref}
\usepackage[hyperref,dvipsnames]{xcolor}
\usepackage[T1]{fontenc} 							
\usepackage[english]{babel} 	
\usepackage{amssymb,amsmath}
\usepackage{cases} 					
\usepackage{bm} 					
\usepackage{stmaryrd} 				
\usepackage{thmbox}
\usepackage{shadethm}		
\usepackage{scalerel} 				
\usepackage{nicefrac}								
\usepackage{siunitx}  		
\usepackage{booktabs} 				
\usepackage{subfigure}		
\usepackage[font={small}]{caption}  
\usepackage{mathtools}				
\usepackage{empheq}					
\usepackage{varwidth}  
\usepackage[normalem]{ulem}

\usepackage{lipsum}

\usepackage{titlesec}				

\titleformat{\section}{\fontsize{10.5}{17}\bfseries}{\thesection}{1em}{}
\titleformat{\subsection}{\bfseries\itshape}{\thesubsection}{1em}{}		

\addto\captionsenglish{}

\DeclareMathAlphabet{\mathcal}{OMS}{cmsy}{m}{n} 

\newcommand\redsout{\bgroup\markoverwith{\textcolor{red}{\rule[0.5ex]{2pt}{0.5pt}}}\ULon}

\newlength\bshft
	\def\fakebold#1{\ThisStyle{\ooalign{$\SavedStyle#1$\cr%
  	\kern-\bshft$\SavedStyle#1$\cr%
  	\kern\bshft$\SavedStyle#1$}}}

\newcommand{\R}{\mathbb{R}}

\newcommand\Pk[2]{{ \mathbb{P}_{#1}{#2} }}

\newcommand\QQdk[2]{{ \fakebold{\mathbb{Q}}_{#1}^\mathrm{dc}{#2} }}

\newcommand\RTk[2]{{ \fakebold{\mathbb{RT}}_{\sqb{#1}}{#2} }}

\newcommand{\dvg}{{ \mathrm{div} }}

\newcommand{\OMEGA}{{ \rb{\Omega} }}

\newcommand\Rey{\mbox{\textit{Re}}}  

\makeatletter
\newcommand*\bigcdot{\mathpalette\bigcdot@{.635}}
\newcommand*\bigcdot@[2]{\mathbin{\vcenter{\hbox{\scalebox{#2}{$\m@th#1\bullet$}}}}}
\makeatother

\DeclareMathOperator{\ip}{{\boldsymbol{\cdot}}}
\DeclareMathOperator{\Fip}{{\boldsymbol{:}}}
\DeclareMathOperator{\DIV}{\nabla\ip}

\DeclareMathOperator{\spn}{{ span }}

\newcommand{\uu}{{ \boldsymbol{u} }}

\newcommand{\vv}{{ \boldsymbol{v} }}

\newcommand{\uuhat}{{ \boldsymbol{\widehat{u}} }}

\newcommand{\xx}{{ \boldsymbol{x} }}
\newcommand{\cc}{{ \boldsymbol{c} }}
\newcommand{\cct}{{ \boldsymbol{\tilde{c}} }}

\newcommand{\sig}{{ \boldsymbol{\sigma} }}
\newcommand{\TAU}{{ \boldsymbol{\tau} }}

\newcommand{\nn}{{ \boldsymbol{n} }}

\newcommand{\ee}{{ \mathbf{e} }}
\newcommand{\zero}{{ \boldsymbol{0} }}
\newcommand{\tend}{{ T }}
\newcommand{\Amat}{{ \boldsymbol{A} }}
\newcommand{\Bmat}{{ \boldsymbol{B} }}
\newcommand{\Dmat}{{ \boldsymbol{D} }}

\newcommand{\drm}{{ \mathrm{d} }}

\newcommand{\mx}{{ \mathrm{max} }}
\newcommand{\phy}{{ \mathrm{phy} }}

\newcommand{\Ctr}[1]{ C_{\mathrm{tr},{#1}} }
\newcommand{\tot}{{ \mathrm{tot} }}
\newcommand{\nm}{{ \mathrm{num} }}
\newcommand{\dxx}{{ \,\drm\xx }}
\newcommand{\dx}{{ \,\drm x }}

\newcommand{\ds}{{\,\drm\boldsymbol{s}}}

\newcommand{\half}{{ \nicefrac{1}{2} }}
\newcommand{\eps}{{ \varepsilon }}

\newcommand{\lavg}{{ \big\{\hspace{-0.99ex}\big\{ }}						
\newcommand{\ravg}{{ \big\}\hspace{-0.99ex}\big\} }}		
\newcommand{\ljmp}{ \left\llbracket }	
\newcommand{\rjmp}{ \right\rrbracket }	
\newcommand\jmp[1]{{ \ljmp#1\rjmp }}							
\newcommand\avg[1]{{ \lavg#1\ravg }}

\newcommand{\nhphantom}[1]{\sbox0{#1}\hspace{-0.751\dimexpr\the\wd0 \relax}}

\newcommand{\Kin}{{ \mathcal{K} }}

\newcommand\Lp[2]{{ L^{#1}{#2} }} 
\newcommand\LP[2]{{ \boldsymbol{L}^{#1}{#2} }}

\newcommand\HM[2]{{ \boldsymbol{H}^{#1}{#2} }}

\newcommand\HDIV{{ \boldsymbol{H}{\rb{\dvg}} }}

\newcommand{\VV}{{ \boldsymbol{V} }}

\newcommand{\Q}{{ Q }}

\newcommand{\T}{{ \mathcal{T}_h }} 
\newcommand{\F}{{ \mathcal{F}_h }}

\newcommand\LIFT[1]{{ \fakebold{\mathcal{L}}{#1} }}	

\newcommand\rb[1]{{ \left(#1\right) }}
\newcommand\sqb[1]{{ \left[ #1 \right] }}
\newcommand\rsb[1]{{ \left(#1\right] }}
\newcommand\set[1]{{ \left\{ #1 \right\} }}

\newcommand\abs[1]{{ \left\lvert#1\right\rvert }}
\newcommand\norm[1]{ \left\lVert#1\right\rVert }

\newcommand\nf[2]{{ \nicefrac{#1}{#2} }}

\newcommand{\tripnorm}[1]{{\left\vert\kern-\nulldelimiterspace\left\vert\kern-\nulldelimiterspace\left\vert #1
	\right\vert\kern-\nulldelimiterspace\right\vert\kern-\nulldelimiterspace\right\vert}}

\newcommand\restr[2]{{												
	\left.\kern-\nulldelimiterspace									
	#1
	\vphantom{\big|}
	\right|_{#2}
	}}
	
\newcommand{\goodgap}{%
	\hspace{0.01\subfigtopskip}
	\hspace{0.01\subfigbottomskip}
	}	

\newtheorem[style=S,underline=true,bodystyle=\normalsize\noindent]{thmDef}{\textsc{Definition}}[section]
\newtheorem[style=S,cut=false]{thmCor}[thmDef]{\textsc{Corollary}}
\newtheorem[style=S,cut=false,headstyle=\normalsize\bfseries\boldmath####1~####2]{thmLem}[thmDef]{\textsc{Lemma}}
 
\newshadetheorem{thm}[thmDef]{\textsc{Theorem}}

\newenvironment{thmProof}
                [1][\unskip]
                { \begin{example}[\normalsize \textsc{Proof #1}] \normalsize}
                { $\hfill\blacksquare$ \end{example} }

\newenvironment{thmEx}
                [0]
                { \refstepcounter{thmDef} \begin{example}[\small\textsc{Example} \thesection.\arabic{thmDef}]  \normalsize}
                {  $\hfill\blacktriangle$ \end{example} }                	
                
\newtheorem[style=S,underline=true,bodystyle=\noindent,cut=false]{thmAss}{\small\textsc{Assumption}}

\journal{`Journal of Computational Physics'}

\definecolor{mediumblue}{RGB}{0,0,205}
\definecolor{forestgreen}{RGB}{34,139,34}
\definecolor{darkred}{RGB}{200,0,0}

\begin{document}

\hypersetup{
  linkcolor=darkred,
  urlcolor=forestgreen,
  citecolor=mediumblue
} 

\begin{frontmatter}


\title{A natural decomposition of viscous dissipation in \\ DG methods for turbulent incompressible flows}

\author[goe]{Christoph Lehrenfeld\corref{cor1}}
\ead{lehrenfeld@math.uni-goettingen.de}
\cortext[cor1]{Corresponding author}
\fntext[fn4]{ORCID: \url{https://orcid.org/0000-0003-0170-8468}}

\author[goe]{Gert Lube}
\ead{lube@math.uni-goettingen.de}  

\author[goe]{Philipp W.\ Schroeder}
\ead{p.schroeder@math.uni-goettingen.de}

\address[goe]{Institute for Numerical and Applied Mathematics, Georg-August-Universit\"at G\"ottingen, 37083 G\"ottingen, Germany \vspace{-7mm}}
\fntext[fn3]{ORCID: \url{https://orcid.org/0000-0001-7644-4693}}




\end{frontmatter}


\vspace{-7mm}
\section{Introduction}	
\label{sec:Introduction}

In this note we aim at a characterisation of the discretisation of viscous dissipation which allows to distinguish `physical' (also frequently called `molecular', or `resolved') from `numerical' dissipation in DG-discretised incompressible flow simulations.
~\\

Let us consider an incompressible Navier-Stokes problem without acting outer forcing terms, i.e.\ without any additional volume forces, and periodic or no-stress boundary conditions.
Given a physical domain $\Omega\subset\R^d$, the strong form of such a problem, equipped with a suitable initial condition $\uu_0\colon \Omega\to\R^d$, reads 
\begin{equation}\label{eq:TINS}
	\partial_t\uu - \nu \Delta \uu + \rb{\uu \ip \nabla} \uu + \nabla p  =  \zero, \qquad 
	\DIV \uu  = 0.
\end{equation}

Here, $\uu \colon\rb{0,\tend}\times\Omega\to\R^d$ indicates the velocity field, $p\colon\rb{0,\tend}\times\Omega\to\R$ is the (zero-mean) kinematic pressure, and the underlying fluid is assumed to be Newtonian with kinematic viscosity $0<\nu\ll 1$.
We are especially interested in the situation where the corresponding Reynolds number is large enough such that a turbulent flow is expected and its approximation is performed in a strongly under-resolved setting.
~\\

For a finite element (FE) pair $\VV_h/Q_h$ for velocity/pressure, and assuming that the simulation is performed up to the time instance $\tend>0$, a typical DG scheme (in primal form) for discretising \eqref{eq:TINS} is written as follows:
\begin{subequations} \label{eq:DG-TINS}
	\begin{empheq}[left=\empheqlbrace]{align} 
	&\text{Find }\rb{\uu_h,p_h}\colon \rsb{0,\tend} \to \VV_h \times \Q_h
	\text{ with }\uu_h\rb{0} =\uu_{0h}\text{ s.t., }\forall\,\rb{\vv_h,q_h}\in\VV_h \times \Q_h,\\
		&\rb{\partial_t\uu_h,\vv_h} 
		+ \nu a_h\rb{\uu_h,\vv_h}
		+ c_h\rb{\uu_h;\uu_h,\vv_h}
		+ j_h\rb{\uu_h,\vv_h}		
		+ b_h\rb{\vv_h,p_h}
		- b_h\rb{\uu_h,q_h}
		= 0.		
	\end{empheq} 
\end{subequations}

The bilinear form $a_h$ treats the viscosity effects, $c_h$ the nonlinear convection term, $b_h$ connects pressure and incompressibility condition and $j_h$ is a possible additional stabilisation and/or turbulence model \cite{Riviere08,PietroErn12}. 
In this note we will focus on the viscous term $a_h$.

\vspace{-5mm}
\section{Physical and numerical dissipation}	
\label{sec:Dissipation}

Testing \eqref{eq:DG-TINS} symmetrically with $\rb{\vv_h,q_h} = \rb{\uu_h,p_h}$ leads to the discrete kinetic energy balance
\begin{equation}
	- \partial_t \Kin\rb{\uu_h}
		= - \frac{\drm}{\drm t} \frac{1}{2}  \norm{\uu_h}_\LP{2}{\OMEGA}^2 
		= \nu a_h\rb{\uu_h,\uu_h} 
		+ c_h\rb{\uu_h;\uu_h,\uu_h} 
		+ j_h\rb{\uu_h,\uu_h},
\end{equation}

whereas on the continuous level, the counterpart for the exact solution $\uu$ is
\begin{equation}
	- \partial_t \Kin\rb{\uu}
		= - \frac{1}{2} \frac{\drm}{\drm t}  \norm{\uu}_\LP{2}{\OMEGA}^2 
		= \nu \norm{\nabla\uu}_\LP{2}{\OMEGA}^2.
\end{equation}

Hence, the only physical dissipation process present in the original Navier--Stokes model is due to viscosity.
Therefore, in our opinion, every additional energy-dissipating (or even energy-producing) mechanism, which is frequently incorporated in $c_h$ and $j_h$, has to be characterised as an artificial (numerical) contribution.
~\\

The purpose of this contribution can be explained compactly as follows:
\emph{We want to distinguish physical and numerical viscous dissipation in $a_h$ and aim for an additive decomposition $a_h = a_h^\phy + a_h^\nm$ where both quantities $a_h^\phy$ and $a_h^\nm$ are non-negative (possibly zero) in order to justify the term `dissipation'. }
In the DG literature, denoting by $\nabla_h$ the broken, i.e.\ element-wise gradient, the choice $a_h^\phy\rb{\uu_h,\uu_h}=\norm{\nabla_h\uu_h}_\LP{2}{}^2$ can be found most frequently (see, e.g.\ \cite{WiartEtAl13}).
We will demonstrate that this definition can be misleading when an under-resolved simulation is performed, and propose an alternative for a large class of DG methods.
~\\

To introduce a mathematically rigorous notion of viscous dissipation processes, let $a_h^\phy$ denote the non-negative part in the discretisation of the viscous term that represents \emph{physical dissipation}, which is supposed to fulfil $a_h^\phy\rb{\uu,\uu}=\norm{\nabla \uu}_\LP{2}{}^2$ for the exact solution $\uu$. 
We assume that the remainder of $a_h$ is a non-negative bilinear form $a_h^\nm$ which describes \emph{numerical dissipation} in the discretisation of the viscous term and require that the decomposition is consistent in the sense that $a_h^\nm\rb{\uu_h,\uu_h}$ vanishes for $h/k \to 0$, with $\uu_h$ being a discrete solution converging to $\uu$ as $h/k \to 0$.
Here, $h$ denotes the underlying mesh size and $k$ is the polynomial order of discrete velocities belonging to $\VV_h$.
~\\

Let us emphasise that the requirement that both parts of the decomposition be non-negative is a restriction and disallows some choices for $a_h^\phy$ which may seem intuitive at first glance.
Being able to identify the physical dissipation, the total numerical dissipation $\eps_h^\tot$ of the scheme can be defined as
\begin{equation}
	\eps_h^\tot \coloneqq
		-  \partial_t \Kin\rb{\uu_h} 
		- \nu a_h^\phy\rb{\uu_h,\uu_h}.
\end{equation}

This total numerical dissipation $\eps_h^\tot$ then fulfils the reasonable and widely accepted expectation (for a meaningful discretisation) that it is non-negative; that is, $\eps_h^\tot\geqslant 0$.

\begin{thmEx}
In an $\HM{1}{}$-conforming setting (see, e.g.\ \cite{HuismannEtAl18}) the physical (viscous) dissipation is simply the scaled seminorm $a^\phy\rb{\uu_h,\uu_h} = \nu \norm{\nabla \uu_h}_\LP{2}{}^2$ while $a_h^\nm\rb{\uu_h,\uu_h} = 0$.
This definition is unproblematic since the discrete velocity $\uu_h$ is continuous here.
Numerical dissipation in $\HM{1}{}$-conforming schemes is thus only contained in explicitly added terms such as turbulence modelling and/or convection stabilisation which are collected in $j_h$, or are part of $c_h$.
\end{thmEx}

\begin{thmEx}
Let $\F=\set{F}$ denote the set of all facets of the decomposition $\T=\set{K}$ and $\phi$ be any piecewise smooth (scalar-, vector- or matrix-valued) function with traces from within the interior of $K^\pm$ denoted by $\phi^\pm$, respectively.
Then, we define the jump $\jmp{\cdot}_F$ and average $\avg{\cdot}_F$ operator across a facet $F\in\F$ by $\jmp{\phi}_F= \phi^+-\phi^-$ and $\avg{\phi}_F=\frac{1}{2}\rb{\phi^+ + \phi^-}$.
The viscous bilinear form of the \emph{non-symmetric} interior penalty (NIP) method \cite{Riviere08} for a scalar problem reads as
\begin{subequations}
	\begin{align} 
		a_h\rb{u_h,v_h} 
			&\coloneqq	\int_\Omega \nabla_h u_h \ip \nabla_h v_h \dxx
				+\sum_{F\in\F} \frac{\lambda}{h_F} \oint_F \jmp{u_h} \jmp{v_h} \ds \\
			&\qquad -\sum_{F\in\F} \oint_F \rb{\avg{\nabla_h u_h} \ip \nn_F} \jmp{v_h} \ds
				+\sum_{F\in\F} \oint_F \jmp{u_h} \rb{\avg{\nabla_h v_h} \ip \nn_F} \ds,
	\end{align}		
\end{subequations}

where $\lambda>0$ is the NIP stabilisation parameter and $h_F$ denotes a length scale for the facet $F$, as usual in DG methods.
In this case we have a simple decomposition with
	\begin{align*}
		a_h\rb{u_h,u_h} 
			&=	\int_\Omega \abs{\nabla_h u_h}^2 \dxx
				+\sum_{F\in\F} \frac{\lambda}{h_F} \oint_F \jmp{u_h}^2 \ds
			= a_h^\phy\rb{u_h,u_h} + a_h^\nm\rb{u_h,u_h}.
	\end{align*}	
\end{thmEx}

In contrast to the decomposition in the previous examples, the decomposition in most other DG schemes is more involved. 
Especially, $a_h^\phy\rb{\uu_h,\uu_h} = \norm{\nabla_h \uu_h}_\LP{2}{}^2$ is frequently not a valid option as the remainder of $a_h$ is not necessarily non-negative, as the following example shows.

\begin{thmEx} \label{ex:1D-SIP}
Let us consider the scalar 1D example (with $\nu=1$) where the domain $\Omega = (0,h)$ with $h=1$ is only one element with periodic boundary conditions, and use the symmetric interior penalty (SIP) DG discretisation \cite{PietroErn12} with scalar-valued polynomial space $V_h$ of order $k=1$.
The set of facets is only $\F=\{1\}$ (due to periodicity).
Using $ \jmp{v_h} = v_h(1) - v_h(0)$ and $\avg{v_h} = \frac12 v_h(0) + \frac12 v_h(1)$, the symmetrically tested bilinear form in this case is
	\begin{align*}
		a_h\rb{u_h,u_h} 
			= \int_0^1 (u_h'(x))^2 \dx 
			- 2  \jmp{u_h}(1) \avg{u_h'}(1) 
			+ \lambda \jmp{u_h}^2 (1). 
	\end{align*}
	
Here, $\lambda>0$ is the SIP penalty parameter which needs to be sufficiently large (depending on the constant of a discrete trace inequality) such that $a_h$ defines an inner product on $V_h\backslash{}\R$ and discrete coercivity is ensured.
Choosing $\lambda = \nf{3}{2}>\nf{\Ctr{0}^2}{2} = 1$ is sufficient as shown in \ref{sec:AppendixB}.
Taking $u_h^\ast = x$, we obtain $u_h'=1$, $\avg{u_h'}(1) = 1$ and $\jmp{u_h}(1) = 1$ which results in
$a_h\rb{u_h^\ast,u_h^\ast} = \frac{1}{2}$ and $\norm{\nabla_h u_h^\ast}_\Lp{2}{\OMEGA}^2 = 1$.
Now, the choice $a_h^\phy\rb{u_h^\ast,u_h^\ast} = \norm{ \nabla_h u_h^\ast}_\Lp{2}{}^2 = 1$ renders $a_h^\nm\rb{u_h^\ast,u_h^\ast}=-\frac{1}{2}$ negative, which contradicts our intuitive understanding that both physical and numerical dissipation should be non-negative.
\end{thmEx}

We conclude that the choice $a_h^\phy\rb{\uu_h,\uu_h} = \norm{\nabla_h \uu_h}_\LP{2}{}^2$ can be misleading  in DG methods.
However, let us mention that the difference between different notions of physical dissipation in DG methods is only relevant in the under-resolved case.
The remainder of this work will demonstrate that a lifting technique can be used to define a more suitable decomposition of the total viscous dissipation into a physical and a numerical contribution. 
In doing so, we restrict ourselves to the SIP method as a very frequently used DG method.

\vspace{-5mm}
\section{A natural decomposition of viscous dissipation for DG methods}	
\label{sec:NaturalDecomposition}

The SIP bilinear form is given by \cite{PietroErn12}
\begin{subequations}\label{eq:SIP-DG-Flow}
\begin{align}
	a_h\rb{\uu_h,\vv_h} 
		&\coloneqq	\int_\Omega \nabla_h \uu_h \Fip \nabla_h \vv_h \dxx
			+\sum_{F\in\F} \frac{\lambda}{h_F} \oint_F \jmp{\uu_h} \ip \jmp{\vv_h} \ds \\
		&\qquad -\sum_{F\in\F} \oint_F \avg{\nabla_h \uu_h} \nn_F \ip \jmp{\vv_h} \ds
			-\sum_{F\in\F} \oint_F \jmp{\uu_h} \ip \avg{\nabla_h \vv_h} \nn_F \ds,	
\end{align}	
\end{subequations}

where $\lambda>0$ is a sufficiently large (due to a discrete inverse inequality) penalty parameter. 
~\\

We can interpret the DG formulation \eqref{eq:SIP-DG-Flow} in a mixed setting,  cf.\ \cite{ArnoldEtAl02}, which gives a natural definition of a discrete diffusive flux (scaled with $\nu^{-1}$) $\sig_h = \sig_h\rb{\uu_h}$, which is defined element-wise for all $\TAU_h \in \restr{\nabla_h\VV_h}{K}$ by the following operation on any $K\in\T$, cf. \cite[eqn.\ (1.2)]{ArnoldEtAl02}: 
\begin{equation}\label{eq:lift}
	\int_K \sig_h \Fip \TAU_h \dxx 
		= - \int_K \uu_h \ip \rb{\DIV\TAU_h} \dxx
			+ \oint_{\partial K} \uuhat_h \ip \rb{\TAU_h\nn_K} \ds 
		= \int_K \nabla \uu_h \Fip \TAU_h \dxx 
			+ \oint_{\partial K} \rb{\uuhat_h - \uu_h} \ip \rb{ \TAU_h \nn_K} \ds 
\end{equation}	

Here, the second equality is due to integration by parts and $\uuhat_h$ denotes a `numerical trace' which characterises different DG methods, see \cite[Table 3.1]{ArnoldEtAl02}.
In the following, for the sake of simplicity, we exclusively want to consider the SIP method where $\uuhat_h=\avg{\uu_h}$.
For SIP, $\uuhat_h - \uu_h=\avg{\uu_h}-\uu_h=-\nf12\jmp{\uu_h}$ and thus we define the lifting operator $\LIFT{}\colon\restr{\VV_h}{\partial K}\to\restr{\nabla_h\VV_h}{K}$ by
\begin{equation}\label{eq:LiftingSIP}
		\oint_{\partial K} \rb{\uuhat_h - \uu_h} \ip \rb{ \TAU_h \nn_K} \ds
		= - \oint_{\partial K} \jmp{\uu_h} \ip \frac{1}{2}\rb{ \TAU_h \nn_K} \ds
		= - \int_K \LIFT{\rb{\jmp{\uu_h}}}\Fip \TAU_h \dxx,
		\quad \forall\,K\in\T.
\end{equation}

With the notion \eqref{eq:LiftingSIP} of the lifting operator $\LIFT{}$, \eqref{eq:lift} can finally be used to obtain the characterisation
\begin{equation} \label{eq:sigmah}
	\sig_h\rb{\uu_h} = \nabla_h \uu_h -  \LIFT{\rb{\jmp{\uu_h}}}.
\end{equation}
 
With this definition of $\sig_h$, one can rewrite the symmetrically tested bilinear form $a_h$ from \eqref{eq:SIP-DG-Flow} as follows:
\begin{equation} \label{eq:NewDecomposition}
	a_h\rb{\uu_h,\uu_h} 
		=\underbrace{ \int_\Omega \abs{\sig_h}^2 \dxx }_{=a_h^\phy\rb{\uu_h,\uu_h}}
			+ \underbrace{\sum_{F\in\F} \frac{\lambda}{h_F} \oint_F \abs{\jmp{\uu_h}}^2 \ds 			
			- \int_\Omega \abs{\LIFT{\rb{\jmp{\uu_h}}}}^2 \dxx}_{=a_h^\nm\rb{\uu_h,\uu_h}}
\end{equation}

We notice that the usual assumption on the parameter $\lambda$ guarantees that both parts $a_h^\phy\rb{\uu_h,\uu_h}$ and $a_h^\nm\rb{\uu_h,\uu_h}$ are non-negative for any discrete function $\uu_h\in\VV_h$; a detailed explanation for this statement can be found in \ref{sec:AppendixB}.
Further, note that $a_h^\nm\rb{\uu_h,\uu_h}\to 0$ as $h/k\to 0$.
~\\

Explaining how \eqref{eq:NewDecomposition} emerges, as shown in detail in \ref{sec:AppendixA}, one can rewrite $a_h^\phy\rb{\uu_h,\uu_h}$ as
\begin{equation} \label{eq:BR}
	a_h^\phy\rb{\uu_h,\uu_h}
		=	\int_\Omega \abs{\nabla_h \uu_h}^2 \dxx
			- 2 \sum_{F\in\F} \oint_F  \avg{\nabla_h \uu_h} \nn_F \ip \jmp{\uu_h}  \ds 
			+ \int_\Omega \abs{\LIFT{\rb{\jmp{\uu_h}}}}^2 \dxx.
\end{equation}

Let us comment on a few topics. 
Firstly, the bilinear form $a_h^\phy$ in \eqref{eq:BR} corresponds to the DG method by Bassi and Rebay \cite{BassiRebay97} and can be seen as a central flux approximation to diffusion/viscosity (to the corresponding first order system).
Secondly, note that for SIP a piecewise constant function will not induce physical dissipation if exclusively the broken gradient is used for the definition of $a_h^\phy$. 
In contrast, using the definition \eqref{eq:NewDecomposition} of $a_h^\phy$ proposed here, also piecewise constant functions induce physical dissipation.
~\\

Moreover, the procedure in the definition of a suitable decomposition of hybrid DG (HDG) methods follows the same reasoning.
In the next section, an $\HDIV$-conforming HDG method will be used for the 3D simulations because of its superior effectivity with respect to computational cost.

\vspace{-5mm}
\section{Numerical demonstration: 3D Taylor--Green vortex problem}	
\label{sec:TGV}

We consider the classical 3D Taylor--Green vortex (TGV) problem which is frequently used to investigate the performance of flow solvers for freely decaying turbulence \cite{TaylorGreen37,Brachet91}. 
In the periodic box $\Omega=\rb{0,2\pi}^3$,
\begin{equation} 
	\uu_0\rb{\xx}
	= \rb{\cos\rb{x_1} \sin\rb{x_2} \sin\rb{x_3},-\sin\rb{ x_1} \cos\rb{ x_2} \sin\rb{x_3},0}^\dag,
\end{equation}

the space-periodic initial condition, is the only driving force.
For the subsequent simulations, the considered Reynolds number is $\Rey=\nu^{-1}=\num{1600}$ and the computations are performed until $\tend=20$.
~\\

The domain is decomposed into $N^3$ cubes and the $\HDIV$-conforming Raviart--Thomas element $\RTk{k}{}$ \cite{BoffiEtAl13} is employed in an exactly divergence-free HDG framework similar to \cite{LehrenfeldSchoeberl16}.
Especially, in order to focus on viscous effects, we do not use any convection stabilisation or additional terms and hence, $c_h\rb{\uu_h;\uu_h,\uu_h}=0$ and $j_h\equiv 0$.
Thus, $-\partial_t\Kin\rb{\uu_h}=\nu a_h\rb{\uu_h,\uu_h}$ and $\eps_h^\tot=\nu a_h^\nm\rb{\uu_h,\uu_h}$.
Concerning the SIP penalty term, we are interested in the smallest penalty that guarantees non-negative total dissipation. 
For the (vector-valued) heat equation this constant can be computed explicitly in the case of a periodic Cartesian mesh as $\lambda^*=\Ctr{k}=\rb{k+1}\rb{k+2}$ for the considered $\RTk{k}{}$-HDG method, or $\lambda^*=k \rb{k+1} / 2$ for a $\QQdk{k}{}$-DG method. 
The corresponding calculations can be found in \ref{sec:AppendixB}. 
Note that due to the incompressibility constraint, which results in the fact that $a_h$ only acts on the (discretely) divergence-free subspace of the DG and HDG methods, the actual minimal penalty parameter can be smaller.
A reference solution with $k=8$ and $N=16$ has been computed and we focus on the comparison of viscous dissipation in the under-resolved situation $k=4$, $N=8$ for different SIP penalties $\lambda\in\set{2,1.5,1.25,1}\lambda^*$.
~\\

\begin{figure}[t]
 	\centering
 		\includegraphics[width=0.475\textwidth]
 			{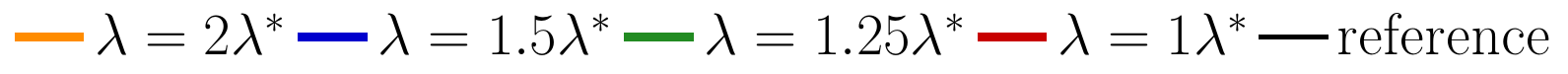} \\ 
 		\includegraphics[width=0.475\textwidth]
			{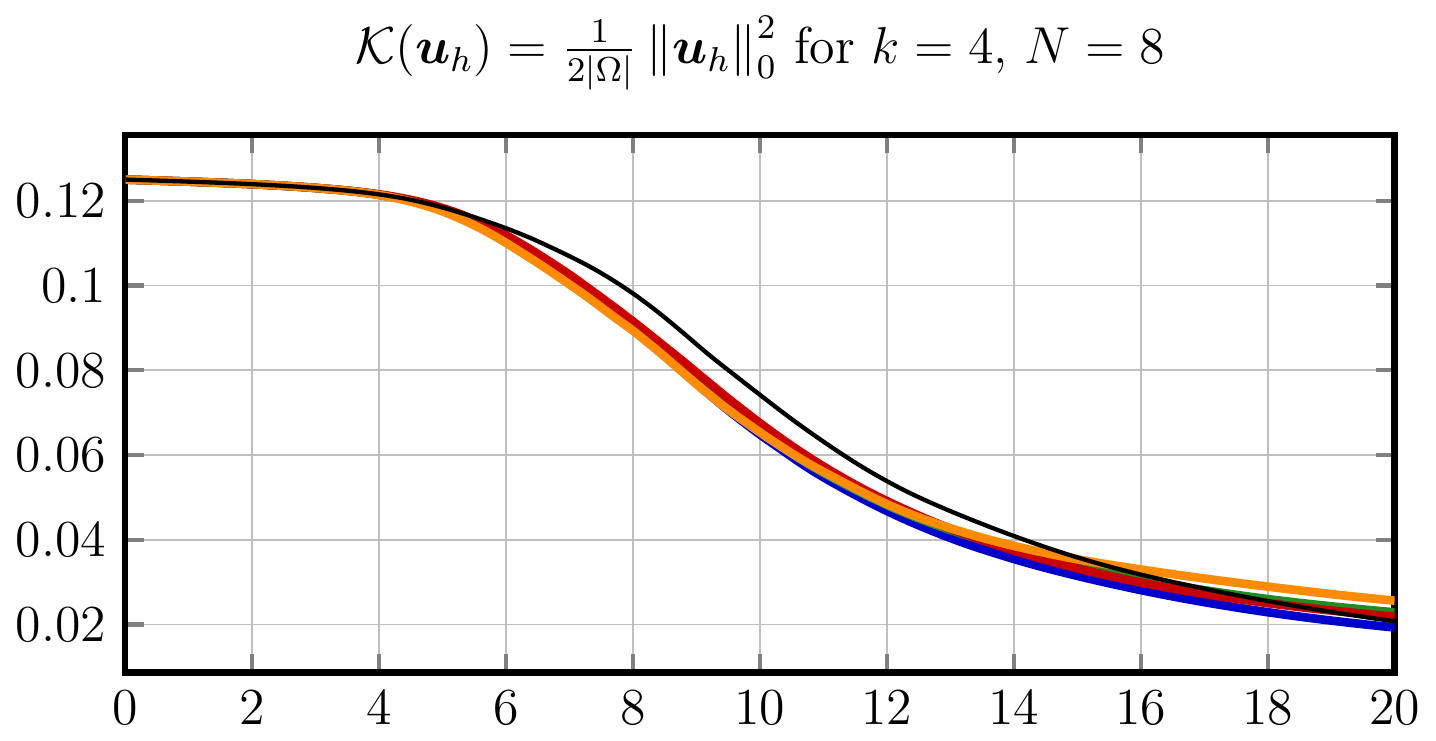} \goodgap
 		\includegraphics[width=0.475\textwidth]
			{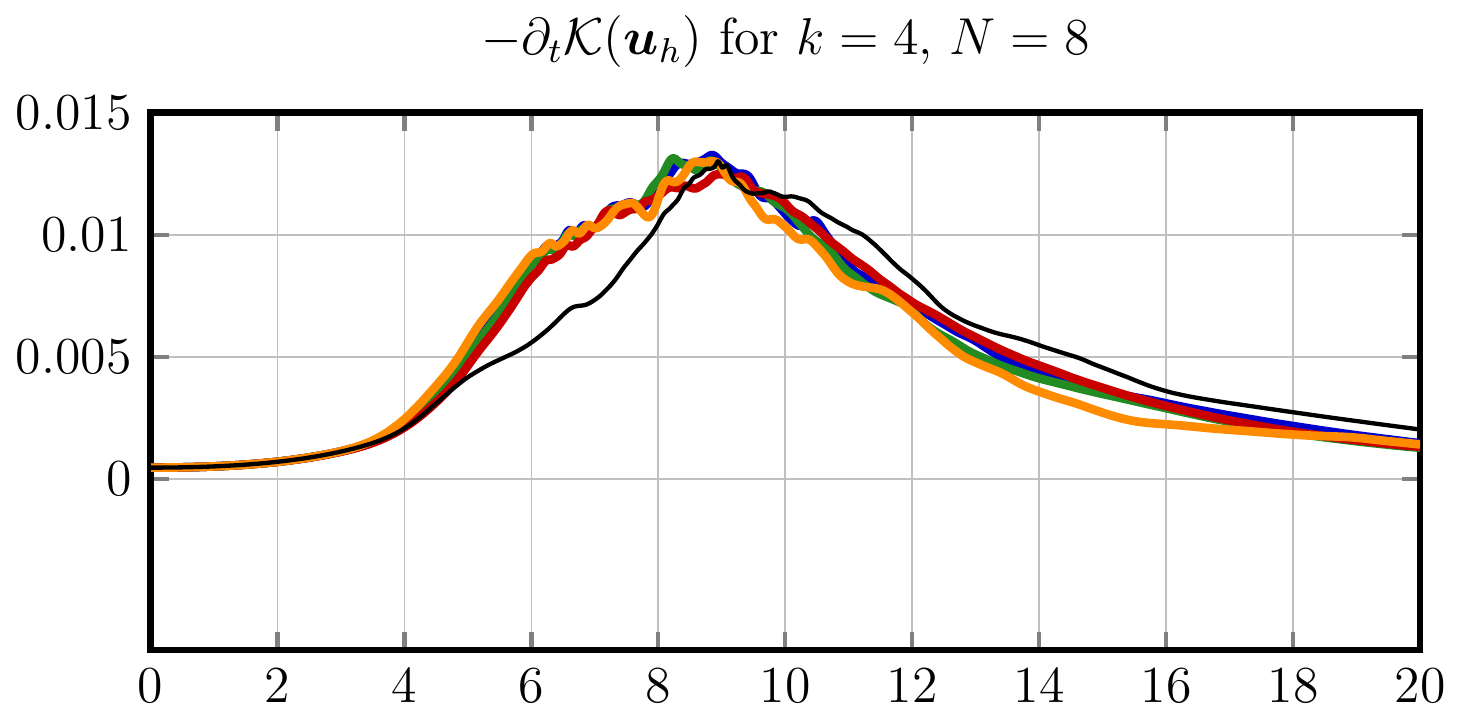} \\
 		\includegraphics[width=0.475\textwidth]
			{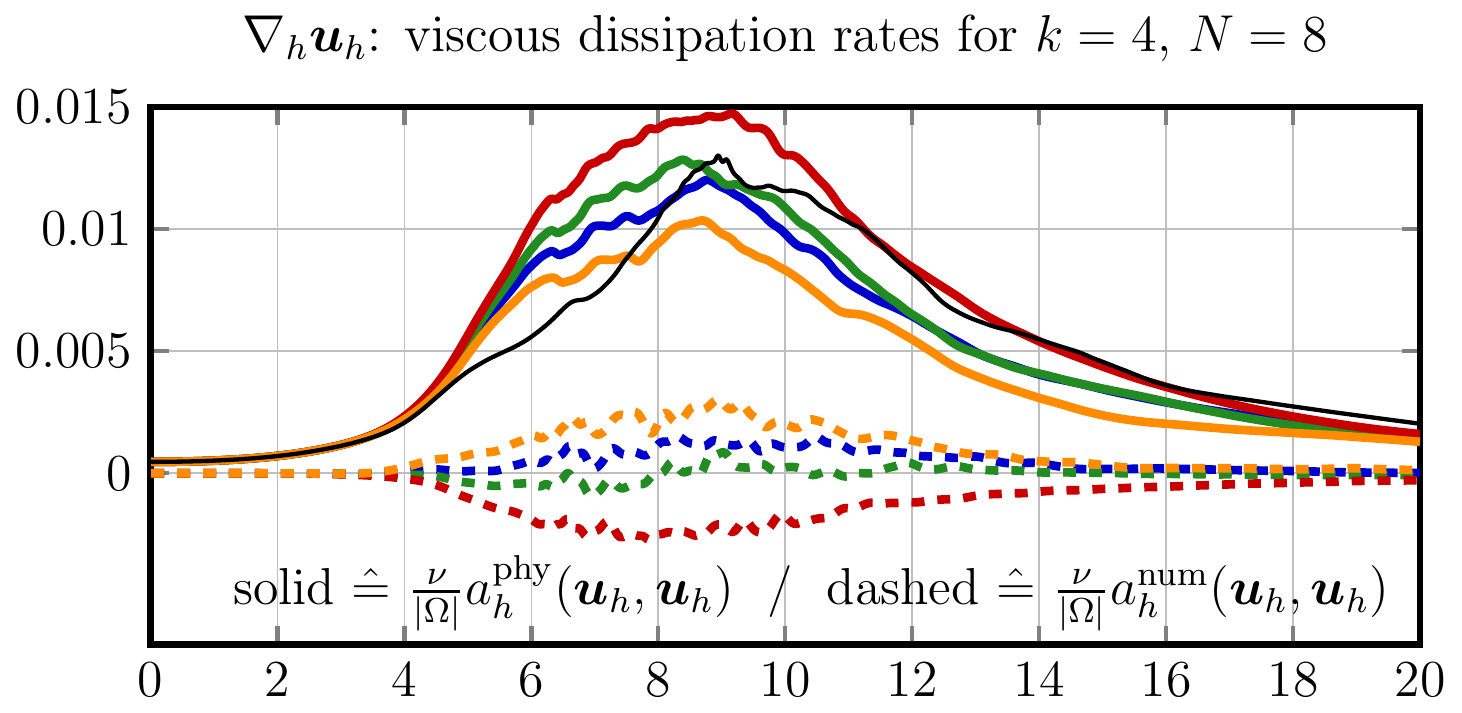} \goodgap
 		\includegraphics[width=0.475\textwidth]
			{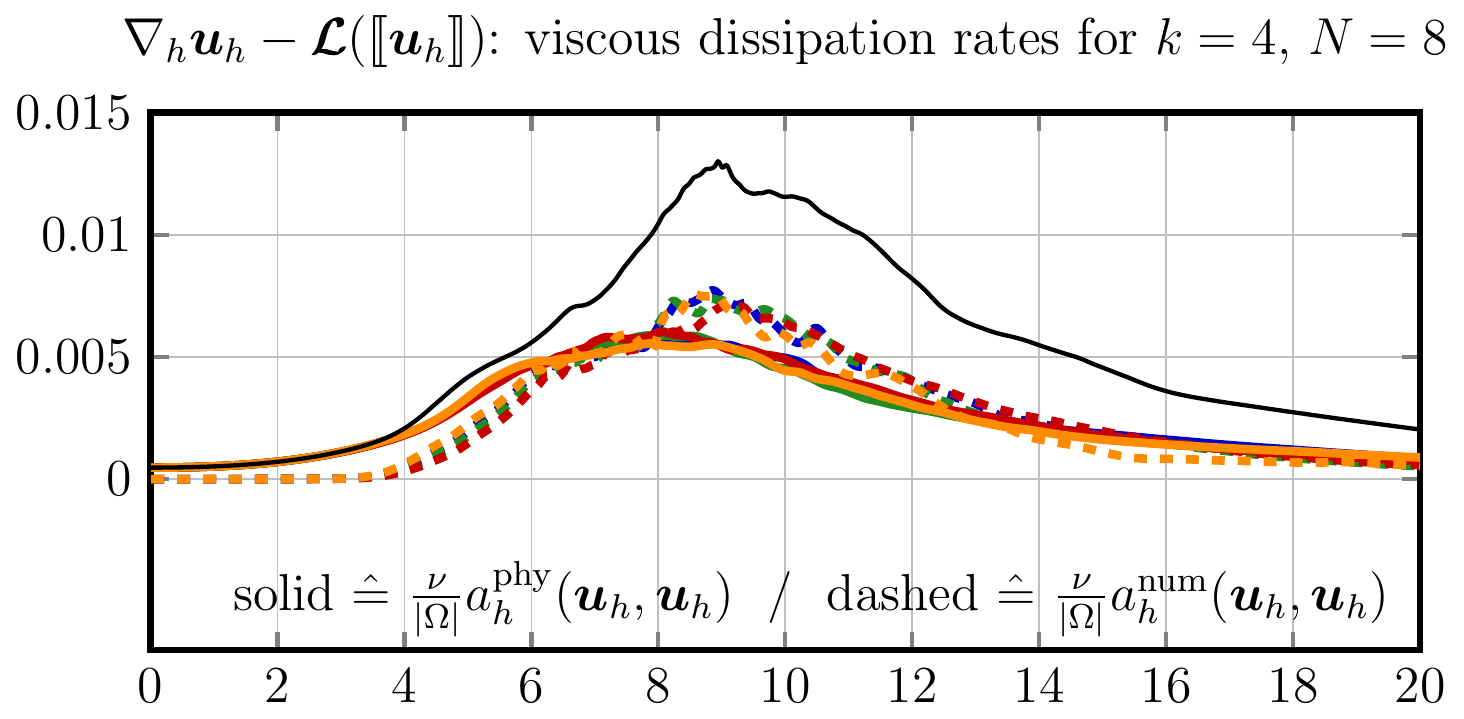}		
 	\caption{Kinetic energy, kinetic energy dissipation rate and viscous dissipations (physical and numerical) for $t\in\sqb{0,20}$, computed with $k=4$, $N=8$. In the bottom row, using the broken gradient $\nabla_h\uu_h$ in the definition of the physical dissipation (left) leads to negative numerical viscous dissipation while using $\sig_h$ from \eqref{eq:sigmah} (right) does not. }
 	\label{fig:3DTGV-EnergyAndDissRates}
 \end{figure}

The upper row of Fig.~\ref{fig:3DTGV-EnergyAndDissRates} shows that, largely unimpaired by the penalty parameter, the evolution of both the kinetic energy $\Kin\rb{\uu_h}$ and the (negative) total kinetic energy dissipation rate $-\partial_t\Kin\rb{\uu_h}$ is reasonable although we are strongly under-resolved.

Concerning the viscous dissipation rates, the bottom row shows how the \emph{interpretation} of physical and numerical dissipation differs if broken gradients are used (left) or if our proposed decomposition \eqref{eq:NewDecomposition} is taken into account (right).
Let us stress that the discretisation is the same for both columns and that for $t\leqslant 4$ the flow seems to be resolved, which renders the difference between which metric is used irrelevant in a resolved simulation.
~\\

For larger times, however, one can observe that with broken gradients the perception of the physical dissipation rate strongly depends on the penalty parameter $\lambda$ even though the total kinetic energy dissipation rate $-\partial_t\Kin\rb{\uu_h}$ does not change much.
Furthermore, if $\lambda$ is chosen small, but still sufficiently large to guarantee non-increasing kinetic energy, the metric `numerical dissipation' defined by broken gradients (left) can become negative. 
This suggests that this characterisation acts unreasonably. 
However, in our opinion this is not a flaw in the method but rather in the metric: There is no conclusion that can be drawn from the sign of this metric.
On the other hand, when our proposed method \eqref{eq:NewDecomposition} of decomposing the viscous dissipation is used (right), both physical and numerical dissipation are non-negative. 
~\\

Lastly, concerning the interpretation of the amount of resolution a simulation can offer, the definition with the broken gradients suggests that most of the total dissipation stems from physical dissipation, whereas our proposition \eqref{eq:NewDecomposition} distributes it more or less evenly between physical and numerical viscous dissipation. 
In view of a clearly under-resolved simulation for $t>4$, the latter behaviour is much more natural to us.

\vspace{-5mm}
\def\bibsection{\section*{References}}
\bibliography{DG-dissipation-BibTeX}
\bibliographystyle{abbrvurl}

\appendix
\vspace{-5mm}
\section{Characterisation of $a_h^\phy$ for SIP}
\label{sec:AppendixA}

Inserting \eqref{eq:sigmah} into \eqref{eq:NewDecomposition}, one directly obtains
\begin{subequations}
\begin{align}
	a_h^\phy\rb{\uu_h,\uu_h}
		&= \int_\Omega \abs{\sig_h}^2 \dxx
		= \rb{\sig_h,\sig_h}_\LP{2}{}
		= \rb{ \nabla_h \uu_h -  \LIFT{\rb{\jmp{\uu_h}}}, \nabla_h \uu_h -  \LIFT{\rb{\jmp{\uu_h}}}}_\LP{2}{} \\
		&= \int_\Omega \abs{\nabla_h \uu_h}^2 \dxx
		- 2 \int_\Omega  \nabla_h \uu_h \Fip \LIFT{\rb{\jmp{\uu_h}}} \dxx
		+ \int_\Omega \abs{\LIFT{\rb{\jmp{\uu_h}}}}^2 \dxx.
\end{align}	
\end{subequations}

The remaining work comes down to a reformulation of the middle term.
In order to accomplish this, the following standard DG procedure can be performed after using definition \eqref{eq:LiftingSIP} of the lifting $\LIFT{\rb{\jmp{\uu_h}}}$:
\begin{subequations} \label{eq:LiftinElBndSkeleton}
\begin{align}
	 \int_\Omega  \nabla_h \uu_h \Fip \LIFT{\rb{\jmp{\uu_h}}} \dxx
		&= \sum_{K\in\T} \oint_{\partial K} \frac{1}{2} \rb{ \nabla_h \uu_h}\nn_K \ip \jmp{\uu_h} \ds \\
		&= \sum_{F\in\F} \oint_{F} 
			\frac{1}{2} \restr{\rb{\nabla_h \uu_h}\nn}{K_1} \ip \jmp{\uu_h} 
			+ \frac{1}{2} \restr{\rb{ \nabla_h \uu_h}\nn}{K_2} \ip \rb{-\jmp{\uu_h}} \ds \\
		&= \sum_{F\in\F} \oint_{F}
			\frac{1}{2} \sqb{\restr{\rb{\nabla_h \uu_h}}{K_1} + \restr{\rb{\nabla_h \uu_h}}{K_2}}\nn_F  \ip \jmp{\uu_h} \ds \\
		&= \sum_{F\in\F} \oint_F  \avg{\nabla_h \uu_h} \nn_F \ip \jmp{\uu_h}  \ds
\end{align}	
\end{subequations}

Note the minus sign in front of the jump term in the second line which stems from the transition of boundary element integrals to the skeleton formulation.

\vspace{-5mm}
\section{Minimal SIP-DG penalty parameter on hyperrectangles}
\label{sec:AppendixB}

In this section, we provide the theoretical foundation for the choice of the SIP penalty parameter in Sec.~\ref{sec:TGV}.

\subsection{Discrete inverse trace inequality and application}
\label{sec:Appendix-TraceInequality}

Firstly, we want to derive a special discrete inverse trace inequality in 1D which takes into account \emph{both} end points of the considered interval.
Such an estimate is crucial in determining a sharp SIP penalty parameter on hyperrectangles.

\begin{thmLem} \label{lem:1DTraceIneq}%
Let $I=\sqb{a,b}$ be an interval with $h=\abs{b-a}$ and $q\in\Pk{k}{\rb{I}}$ be a $k$-th order 1D polynomial.
Then,
\begin{equation}\label{eq:1DInvTrace}
	\abs{q\rb{a}}^2 + \abs{q\rb{b}}^2
		\leqslant \frac{\Ctr{k}^2}{h} \int_a^b \abs{q\rb{x}}^2 \dx
\end{equation}

holds with $\Ctr{k}^2=\rb{k+1}\rb{k+2}$.
\end{thmLem}
\begin{thmProof}
In order to prove the claim, consider the shifted Legendre polynomials $L_m\rb{x}=\widetilde{L}_m\rb{2x-1}$ on the unit interval $\sqb{0,1}$, where $\widetilde{L}_m$ denotes the standard Legendre polynomials defined on $\sqb{-1,1}$.
The crucial properties of these polynomials are $L_m\rb{0}=\rb{-1}^m$, $L_m\rb{1}=1$ and 
\begin{equation}
	\int_0^1 L_m\rb{x} L_n\rb{x} \dx = \frac{1}{2m+1} \delta_{mn}.	
\end{equation}

Thus, using the basis representation $q\rb{x}=\sum_{m=0}^k c_m L_m\rb{x}$, \eqref{eq:1DInvTrace} can be rewritten as
\begin{equation}
	\sum_{m,n=0}^k c_m c_n \sqb{ \rb{-1}^{m+n} + 1 }
		\leqslant \Ctr{k}^2 \sum_{m=0}^k c_m^2 \frac{1}{2m+1}.
\end{equation} 

Defining $\cc=\rb{c_0,\dots,c_k}^\dag$, the matrix $\Bmat=\rb{B_{mn}}_{m,n=0}^k$ with $B_{mn}=\rb{-1}^{m+n} + 1$, and the matrix $\Dmat=\mathrm{diag}\rb{\nf{1}{\rb{2m+1}}}$, one obtains
\begin{equation}
	\cc^\dag \Bmat \cc \leqslant \Ctr{k}^2 \cc^\dag \Dmat \cc 
		\quad \Leftrightarrow \quad 
		\cct^\dag \Amat \cct \leqslant \Ctr{k}^2 \cct^\dag \cct,
\end{equation}

where $\cct=\Dmat^\half \cc$ and $\Amat = \Dmat^{-\half} \Bmat \Dmat^{-\half}$ with $\Dmat^{-\half}=\mathrm{diag}\rb{\sqrt{2m+1}}$.
The entries $A_{mn}$ of $\Amat$ are
\begin{equation}
	A_{mn} = 
	\begin{cases}
	0, & \text{if }m+n\text{ odd,} \\
	4\sqrt{m+\nf12}\sqrt{n+\nf12}, & \text{otherwise.}
	\end{cases}
\end{equation}

Thus, by using the concept of the Rayleigh quotient, determining $\Ctr{k}^2$ reduces to finding the maximum eigenvalue $\lambda_\mx\rb{\Amat}$ of $\Amat$.
More precisely, one can verify that
\begin{equation}
	\Ctr{k}^2 \geqslant \lambda_\mx\rb{\Amat} = \rb{k+1}\rb{k+2}.
\end{equation}

Applying a standard scaling argument and the width $h$, the result can be transferred from $\sqb{0,1}$ to $\sqb{a,b}$.
Furthermore, $\Ctr{k}^2=\rb{k+1}\rb{k+2}$ is the smallest possible constant for which \eqref{eq:1DInvTrace} holds.
\end{thmProof}

Note that we verified Lem.~\ref{lem:1DTraceIneq} numerically and observed that the given $\Ctr{k}^2$ is indeed sharp.
~\\

The second aim is to apply Lem.~\ref{lem:1DTraceIneq} in the special situation where the normal gradient of the velocity, or the lifting operator, in normal direction on facets has to be estimated.
This is the typical application of the discrete trace inequality in DG methods for diffusive problems.
As we want to especially treat the Raviart--Thomas $\RTk{k}{}$ case, we exploit the inclusion $\QQdk{k}{}\subset\RTk{k}{}\subset\QQdk{k+1}{}$ \cite{BoffiEtAl13}.

\begin{thmLem} \label{lem:TraceIneq-NormalGradient}%
For all $\vv_h\in\QQdk{k+1}{}$ and with $\Ctr{k}^2=\rb{k+1}\rb{k+2}$, the following discrete trace inequality holds:
\begin{equation}
	\norm{\rb{\nabla\vv_h}\nn_K}_\LP{2}{\rb{\partial K}}^2
		\leqslant \frac{\Ctr{k}^2}{h_K} \norm{\nabla \vv_h}_\LP{2}{\rb{K}}^2,
		\quad \forall\,K\in\T
\end{equation}
\end{thmLem}
\begin{thmProof}
This proof is performed for 3D; for 2D and 1D, the same result holds and can be shown as a simplification of the 3D case.
Let $K$ be a cube with $K=I_1\times I_2\times I_3=\sqb{a_1,b_1}\times\sqb{a_2,b_2}\times\sqb{a_3,b_3}$ and $h_K=\abs{b_1-a_1}=\abs{b_2-a_2}=\abs{b_3-a_3}$.
The boundary of the cube can be decomposed by means of $\partial K=\bigcup_{i=1}^3 \partial K_i$ with $\partial K_i=\set{\xx\in\partial K\colon \nn_K\rb{\xx} \parallel \ee_i}$, where $\xx=\rb{x_1,x_2,x_3}^\dag$  and $\ee_i$ denotes the Euclidean unit vector in direction $i$.
Then, one obtains
\begin{equation} \label{eq:Qkp1Trace-step0}
	\norm{\rb{\nabla\vv_h}\nn_K}_\LP{2}{\rb{\partial K}}^2
		= \sum_{i=1}^3 \norm{\partial_{x_i}\vv_h}_\LP{2}{\rb{\partial K_i}}^2.
\end{equation}

Note that for $\vv_h\in\QQdk{k+1}{}$ with $\vv_h=\rb{v_1,v_2,v_3}^\dag$, the charatersiation
\begin{equation} \label{eq:Qkp1-characterisation}
	\restr{\vv_h}{K}
		\in 
	\spn \set{ \begin{pmatrix}
 	\Pk{k+1}{} \otimes \Pk{k+1}{} \otimes \Pk{k+1}{} \\
 	\Pk{k+1}{} \otimes \Pk{k+1}{} \otimes \Pk{k+1}{} \\
 	\Pk{k+1}{} \otimes \Pk{k+1}{} \otimes \Pk{k+1}{} \\
 	\end{pmatrix}
 	}
 	\quad \Rightarrow \quad 
	\restr{\rb{\partial_{x_i} v_m}_{i=1}^3}{K}
		\in 
	\spn \set{ \begin{pmatrix}
 	\Pk{k}{} \otimes \Pk{k+1}{} \otimes \Pk{k+1}{} \\
 	\Pk{k+1}{} \otimes \Pk{k}{} \otimes \Pk{k+1}{} \\
 	\Pk{k+1}{} \otimes \Pk{k+1}{} \otimes \Pk{k}{} \\
 	\end{pmatrix}
 	} 	
\end{equation}

is valid for all $m=1,2,3$.
Now, considering for example $i=1$ leads to
\begin{equation} \label{eq:Qkp1Trace-step1}
	\norm{\partial_{x_1}\vv_h}_\LP{2}{\rb{\partial K_1}}^2
		= \int_{I_2} \int_{I_3}
			\abs{\partial_{x_1}\vv_h\rb{a_1,x_2,x_3}}^2 + \abs{\partial_{x_1}\vv_h\rb{b_1,x_2,x_3}}^2
			\dx_2 \dx_3.
\end{equation}

Fortunately, from \eqref{eq:Qkp1-characterisation} we can infer that $\partial_{x_1}\vv_h\rb{\ip,x_2,x_3}\in\sqb{\Pk{k}{\rb{K}}}^3$ and thus, the 1D result from Lem.~\ref{lem:1DTraceIneq} can be applied componentwise: 
\begin{equation}
	\abs{\partial_{x_1}\vv_h\rb{a_1,x_2,x_3}}^2 + \abs{\partial_{x_1}\vv_h\rb{b_1,x_2,x_3}}^2
		\leqslant \frac{\Ctr{k}^2}{h_K} \int_{I_1}
			\abs{\partial_{x_1}\vv_h\rb{x_1,x_2,x_3}}^2
			\dx_1
\end{equation}

Inserting this estimate into \eqref{eq:Qkp1Trace-step1} leads to
\begin{equation}
	\norm{\partial_{x_1}\vv_h}_\LP{2}{\rb{\partial K_1}}^2
		\leqslant \frac{\Ctr{k}^2}{h_K}
			\int_{I_1} \int_{I_2} \int_{I_3} \abs{\partial_{x_1}\vv_h\rb{x_1,x_2,x_3}}^2 \dx_1 \dx_2 \dx_3
			= \frac{\Ctr{k}^2}{h_K} \norm{\partial_{x_1}\vv_h}_\LP{2}{\rb{K}}^2.
\end{equation}

Finally, using the same arguments also for $i=2,3$ and inserting the particular estimates for $\norm{\partial_{x_i}\vv_h}_\LP{2}{\rb{\partial K_i}}^2$ into \eqref{eq:Qkp1Trace-step0} concludes the proof.
\end{thmProof}

\subsection{Non-negativity of numerical viscous dissipation}
\label{sec:Appendix-NonNegViscDiss}

We now want to show that provided a certain minimum SIP penalty parameter $\lambda^*$ is chosen, the numerical viscous dissipation $a_h^\nm\rb{\vv_h,\vv_h}$, defined in \eqref{eq:NewDecomposition}, is non-negative for all $\vv_h\in\VV_h$.
Here, only $\VV_h=\QQdk{k+1}{}$ is considered which includes the $\RTk{k}{}$ case.
For the sake of brevity, suppose we are working only on meshes containing lines/squares/cubes, then $h_F=h_K$ for all $F\in\F$ and for all $K\in\T$.
~\\

\begin{thmLem} \label{lem:ahnum-nonnegative}%
Provided $\lambda\geqslant\lambda^*=\frac{1}{2}\Ctr{k}^2$, the numerical dissipation of \eqref{eq:SIP-DG-Flow} is non-negative; that is,
\begin{equation}
	a_h^\nm\rb{\vv_h,\vv_h} 
		=\sum_{F\in\F} \frac{\lambda}{h_F} \oint_F \abs{\jmp{\vv_h}}^2 \ds 			
			- \int_\Omega \abs{\LIFT{\rb{\jmp{\vv_h}}}}^2 \dxx
		\geqslant 0, \quad  \forall\,\vv_h\in\VV_h=\QQdk{k+1}{}.
\end{equation}
\end{thmLem}
\begin{thmProof}
Rewriting the penalty term from skeleton to boundary element formulation, one obtains
\begin{equation}
	a_h^\nm\rb{\vv_h,\vv_h}
		= \sum_{K\in\T} \frac12 \frac{\lambda}{h_K} \oint_{\partial K} \abs{\jmp{\vv_h}}^2 \ds
			- \sum_{K\in\T}\int_K \abs{\LIFT{\rb{\jmp{\vv_h}}}}^2 \dxx.
\end{equation}

Inserting definition \eqref{eq:LiftingSIP} for the SIP lifting, and using Cauchy--Schwarz and Young ($\eps>0$), the estimate
\begin{align}
	\sum_{K\in\T}\int_K \abs{\LIFT{\rb{\jmp{\vv_h}}}}^2 \dxx
		&= \sum_{K\in\T} \oint_{\partial K} \jmp{\vv_h} \ip \frac12 \LIFT{\rb{\jmp{\vv_h}}}\nn_K \ds \\
		&\leqslant 
			\sum_{K\in\T} \frac{\eps}{4} \oint_{\partial K} \abs{\jmp{\vv_h}}^2 \ds
			+ \sum_{K\in\T} \frac{1}{4\eps} \oint_{\partial K} \abs{\LIFT{\rb{\jmp{\vv_h}}}\nn_K}^2 \ds
\end{align}

holds.
Furthermore, due to the fact that $\LIFT{}\colon\restr{\VV_h}{\partial K}\to\restr{\nabla_h\VV_h}{K}$, one can now apply the discrete trace inequality (Lem.~\ref{lem:TraceIneq-NormalGradient}) to infer
\begin{equation} 
	\oint_{\partial K} \abs{\LIFT{\rb{\jmp{\vv_h}}}\nn_K}^2 \ds
		\leqslant \frac{\Ctr{k}^2}{h_K} \norm{\LIFT{\rb{\jmp{\vv_h}}}}_\LP{2}{\rb{K}}^2.
\end{equation}

Choosing $\eps={\Ctr{k}^2}/{(2h_K)}$ and reordering thus yields
\begin{equation} \label{eq:LiftingBoundedness}
	\sum_{K\in\T}\int_K \abs{\LIFT{\rb{\jmp{\vv_h}}}}^2 \dxx
		\leqslant 
			\sum_{K\in\T} \frac{\Ctr{k}^2}{4h_K} \oint_{\partial K} \abs{\jmp{\vv_h}}^2 \ds.
\end{equation}

Inserting this estimate into the definition of $a_h^\nm$ leads to
\begin{align}
	a_h^\nm\rb{\vv_h,\vv_h}
		&\geqslant \sum_{K\in\T} \frac12 \frac{\lambda}{h_K} \oint_{\partial K} \abs{\jmp{\vv_h}}^2 \ds
			- \sum_{K\in\T} \frac{\Ctr{k}^2}{4h_K} \oint_{\partial K} \abs{\jmp{\vv_h}}^2 \ds \\
		&\geqslant \sum_{K\in\T}
			\frac12\rb{\frac{\lambda-\frac12\Ctr{k}^2}{h_K}} \oint_{\partial K} \abs{\jmp{\vv_h}}^2 \ds.
\end{align}

Concluding the proof, whenever $\lambda$ is chosen according to $\lambda\geqslant\lambda^*=\frac{1}{2}\Ctr{k}^2$,  $a_h^\nm$ is non-negative.
\end{thmProof}

\subsection{Minimal SIP parameter}
\label{sec:Appendix-MinSIPParameter}

Lastly, we establish a connection between non-negativity of $a_h^\nm$ and the discrete coercivity (stability) of the corresponding SIP-DG method \eqref{eq:SIP-DG-Flow}.

\begin{thmLem} \label{lem:ahnum-coercive}%
Provided $\lambda>\lambda^*=\frac{1}{2}\Ctr{k}^2$, the SIP-DG method \eqref{eq:SIP-DG-Flow} is coercive on $\VV_h=\QQdk{k+1}{}$.
\end{thmLem}
\begin{thmProof}
Testing \eqref{eq:SIP-DG-Flow} symmetrically, inserting the definition of the lifting operator and going over from skeleton to boundary element formulation (using \eqref{eq:LiftinElBndSkeleton}) yields
\begin{subequations}
\begin{align}
	a_h\rb{\vv_h,\vv_h} 
		&=	\int_\Omega \abs{\nabla_h \vv_h}^2 \dxx
			-2\sum_{F\in\F} \oint_F \avg{\nabla_h \vv_h} \nn_F \ip \jmp{\vv_h} \ds
			+\sum_{F\in\F} \frac{\lambda}{h_F} \oint_F \abs{\jmp{\vv_h}}^2 \ds \\
		&= \sum_{K\in\T} \int_K \abs{\nabla \vv_h}^2 \dxx
			-2\sum_{K\in\T} \int_K  \nabla \uu_h \Fip \LIFT{\rb{\jmp{\uu_h}}} \dxx
			+\sum_{K\in\T} \frac12 \frac{\lambda}{h_K} \oint_{\partial K} \abs{\jmp{\vv_h}}^2 \ds.
\end{align}		
\end{subequations}

Applying Cauchy--Schwarz and Young's inequality ($\eps>0$) to the problematic middle term and exploiting the boundedness of the lifting operator \eqref{eq:LiftingBoundedness}, one obtains
\begin{equation}
	\sum_{K\in\T} \int_K  \nabla \uu_h \Fip \LIFT{\rb{\jmp{\uu_h}}} \dxx
		\leqslant 
			\sum_{K\in\T} \frac{1}{2\eps} \int_K \abs{\nabla_h \vv_h}^2 \dxx
			+ \sum_{K\in\T} \frac{\eps}{2} \frac{\Ctr{k}^2}{4h_K} \oint_{\partial K} \abs{\jmp{\vv_h}}^2 \ds.
\end{equation}

Inserting this estimate leads to
\begin{equation}
	a_h\rb{\vv_h,\vv_h}
		\geqslant \rb{1-\frac{1}{\eps}} \norm{\nabla_h\vv_h}_\LP{2}{}^2
		 	+\sum_{K\in\T} \frac12 \rb{ \frac{\lambda-\frac{\eps}{2}\Ctr{k}^2}{h_K} } 
		 		\oint_{\partial K} \abs{\jmp{\vv_h}}^2 \ds.	
\end{equation}

Choosing $\eps>0$ infinitesimal, the minimum stabilisation can be achieved by $\lambda> \lambda^*=\frac{1}{2}\Ctr{k}^2$.
This coincides with the minimum SIP penalty parameter which is needed for the non-negativity of $a_h^\nm$.
\end{thmProof}

\end{document}